\documentclass[10pt]{amsart}
\usepackage{amsmath,amssymb,latexsym,hyperref,url,graphicx}

\newtheorem{theorem}{Theorem}[section]

\numberwithin{equation}{section}
\theoremstyle{definition}
\newtheorem{example}[theorem]{Example}

\newcommand{\Z}{\mathbb{Z}}
\newcommand{\legendre}[2]{{#1 \overwithdelims () #2}}
\newcommand{\nequiv}{\mathrel{\not\equiv}}
\newcommand{\rank}{\operatorname*{rank}}
\newcommand{\colonequal}{\mathrel{\mathop:}=}

\begin{document}

\title{$p$-regularity of the $p$-adic valuation of the Fibonacci sequence}

\author{Luis A. Medina}
\address{Department of Mathematics, University of Puerto Rico, San Juan, PR 00931}
\email{luis.medina17@upr.edu}

\author{Eric Rowland}
\address{
	Laboratoire de combinatoire et d'informatique math\'ematique \\
	Universit\'e du Qu\'ebec \`a Montr\'eal \\
	Montr\'eal, QC H2X 3Y7, Canada
}
\curraddr{
	University of Liege \\
	D\'epartement de Math\'ematiques \\
	4000 Li\`ege, Belgium
}




\begin{abstract}
We show that the $p$-adic valuation of the sequence of Fibonacci numbers is a $p$-regular sequence for every prime $p$.
For $p \neq 2, 5$, we determine that the rank of this sequence is $\alpha(p) + 1$, where $\alpha(m)$ is the restricted period length of the Fibonacci sequence modulo $m$.
\end{abstract}

\maketitle

\section{Introduction}
\label{introduction}

Let $F_n$ be the $n$th Fibonacci number.
The sequence $(F_n)_{n \geq 1}$ is
\[
	1, 1, 2, 3, 5, 8, 13, 21, 34, 55, 89, 144, 233, 377, 610, 987, \dots.
\]
Let $a \bmod m$ denote the least nonnegative integer $b$ such that $a \equiv b \mod m$.
It is easy to see that for a fixed $m \geq 1$ the sequence $(F_n \bmod m)_{n \geq 1}$ is eventually periodic.
Namely, there are only $m^2$ possible pairs of consecutive terms, so some pair must occur more than once.
Since the Fibonacci numbers modulo $m$ satisfy a recurrence of order $2$, the sequences beginning at two different positions with the same initial pair coincide.

The sequence $(F_n \bmod m)_{n \geq 1}$ is not just eventually periodic but in fact periodic.
This is because the $n$th term can be determined from terms $n + 1$ and $n + 2$, so the recurrence can be run backward uniquely as well as forward.
We denote the (minimal) period length of $(F_n \bmod m)_{n \geq 1}$ by $\pi(m)$.
Let $\alpha(m)$ be the smallest value of $n \geq 1$ such that $F_n \equiv 0 \mod m$.

\begin{example}
The sequence $(F_n \bmod 3)_{n \geq 1}$ of Fibonacci numbers modulo $3$ is
\[
	1, 1, 2, 0, 2, 2, 1, 0, 1, 1, 2, 0, 2, 2, 1, 0, \dots.
\]
It is periodic with $\pi(3) = 8$ and $\alpha(3) = 4$.
\end{example}

The quantities $\pi(m)$ and $\alpha(m)$ are called, respectively, the \emph{period length} and the \emph{restricted period length} of the Fibonacci sequence modulo $m$.
It is known that $\alpha(m)$ divides $\pi(m)$~\cite[Theorem~3]{Wall}.
For $m \geq 2$, note that $\alpha(m) > 1$ since $F_1 = 1 \nequiv 0 \mod m$.
Throughout the paper we will use the following result~\cite{Robinson,Fulton--Morris}, where $\legendre{a}{b}$ is the Legendre symbol.

\begin{theorem}\label{mod5}
Let $p$ be a prime.
Then $\alpha(p) \mid p - \legendre{5}{p}$.
\end{theorem}

For $p \equiv 1,4 \mod 5$ this implies $\alpha(p) \mid p - 1$, and for $p \equiv 2,3 \mod 5$ it implies $\alpha(p) \mid p + 1$.

Consider integers $k \geq 2$ and $n \geq 1$. The exponent of the highest power of $k$ that divides $n$ is denoted by $\nu_k(n)$.  For example, $\nu_2(144) = 4$.  If $k = p$ is prime, $\nu_p(n)$ is called the \emph{$p$-adic valuation} of $n$.

\begin{example}
The sequence $\nu_3(F_n)_{n \geq 1}$ of $3$-adic valuations of the Fibonacci numbers is
\[
	0, 0, 0, 1, 0, 0, 0, 1, 0, 0, 0, 2, 0, 0, 0, 1, \dots.
\]
\end{example}

Lengyel~\cite{Lengyel 1995} discovered the structure of $\nu_p(F_n)_{n \geq 1}$ for prime $p$.

\begin{theorem}[Lengyel]\label{Lengyel}
Let $p \neq 2, 5$ be a prime and $n \geq 1$.
Then
\begin{align*}
	\nu_5(F_n) &= \nu_5(n), \\
	\nu_2(F_n) &=
	\begin{cases}
		\nu_2(n) + 2	& \text{if $n \equiv 0 \mod 6$} \\
		1		& \text{if $n \equiv 3 \mod 6$} \\
		0		& \text{if $n \equiv 1,2,4,5 \mod 6$,}
	\end{cases} \\
	\nu_p(F_n) &=
	\begin{cases}
		\nu_p(n) + \nu_p(F_{\alpha(p)})	& \text{if $n \equiv 0 \mod \alpha(p)$} \\
		0				& \text{if $n \nequiv 0 \mod \alpha(p)$.}
	\end{cases}
\end{align*}
\end{theorem}

Initial computations reveal that the constant term $\nu_p(F_{\alpha(p)})$ is equal to $1$ for the first few primes.
The statement $\nu_p(F_{\alpha(p)}) = 1$ is equivalent to $\pi(p^2) \neq \pi(p)$~\cite[Corollary~3.33]{Renault}.
A prime $p$ such that $\nu_p(F_{\alpha(p)}) \neq 1$ is called a \emph{Wall--Sun--Sun prime} (sometimes also called a \emph{Fibonacci--Wieferich prime}).
It is not known if any Wall--Sun--Sun primes exist. In 2007, McIntosh and Roettger~\cite{McIntosh--Roettger} showed that there are no Wall--Sun--Sun primes $p$ with $p < 2 \times 10^{14}$.  Dorais and Klyve~\cite{Dorais--Klyve} extended this bound to $9.7\times 10^{14}$.  The current bound, provided by the PrimeGrid project~\cite{primegrid}, is $2.8\times 10^{16}$.

Valuations of various combinatorial sequences have been studied by a number of researchers~\cite{Amdeberhan--Manna--Moll, Amdeberhan--Medina--Moll, Cohn, Lengyel 1994, Postnikov--Sagan, Sun--Moll 2009, Sun--Moll 2010}.
In this paper we study the $p$-adic valuation of the Fibonacci sequence from the perspective of regular sequences, which we now define.

For $k \geq 2$, the \emph{$k$-kernel} of a sequence $s(n)_{n \geq 0}$ is the set of subsequences
\[
	\ker_k s(n)_{n \geq 0} \colonequal \{s(k^e n + i)_{n \geq 0} \, : \, e \geq 0, \, 0 \leq i \leq k^e - 1\}.
\]
A sequence $s(n)_{n \geq 0}$ is \emph{$k$-regular} if the $\Z$-module $\langle \ker_k s(n)_{n \geq 0} \rangle$ generated by its $k$-kernel is finitely generated.
The \emph{rank} of a $k$-regular sequence is the rank of this $\Z$-module.
The rank of a $k$-regular sequence is analogous to the order of a sequence satisfying a linear recurrence with constant coefficients.
For example, since the Fibonacci sequence satisfies the recurrence $F_{n+2} = F_{n+1} + F_n$ of order $2$, the rank of the $\Z$-module generated by the set $\{(F_{n + i})_{n \geq 1} \, : \, i \geq 0\}$ is also $2$.

Regular sequences were introduced by Allouche and Shallit~\cite{Allouche--Shallit} and are a natural class of sequences for the study of valuations.
For example, consider the sequence $\nu_k(n+1)_{n \geq 0}$ (where we index terms beginning with $0$ to match the definition).

\begin{theorem}\label{valuation rank}
Let $k \geq 2$.
Then $\nu_k(n+1)_{n \geq 0}$ is a $k$-regular sequence of rank $2$.
\end{theorem}

\begin{proof}
Let $B = \{\nu_k(n+1)_{n \geq 0}, \nu_k(k (n+1))_{n \geq 0}\}$.
We show that for each sequence $s(n)_{n \geq 0} \in B$, the subsequence $s(k n + i)_{n \geq 0}$ for each $0 \leq i \leq k - 1$ can be written as a $\Z$-linear combination of elements of $B$.
Indeed, we have
\begin{align*}
	\nu_k(k n + i + 1) &=
	\begin{cases}
		0		& \text{if $0 \leq i \leq k-2$} \\
		\nu_k(k (n+1))	& \text{if $i = k-1$,}
	\end{cases} \\
	\nu_k(k (k n + i + 1)) &=
	\begin{cases}
		-\nu_k(n+1) + \nu_k(k (n+1))	& \text{if $0 \leq i \leq k-2$} \\
		-\nu_k(n+1) + 2 \nu_k(k (n+1))	& \text{if $i = k-1$.}
	\end{cases}
\end{align*}
By an induction argument, this implies that every sequence in the $k$-kernel of $\nu_k(n+1)_{n \geq 0}$ is a $\Z$-linear combination of elements of $B$.
Therefore $\nu_k(n+1)_{n \geq 0}$ is $k$-regular with rank at most $2$.
The two sequences in $B$ are linearly independent, since $-\nu_k(n+1) + \nu_k(k (n+1)) = 1$ for all $n \geq 0$.
Since both sequences belong to the $k$-kernel, the rank is exactly $2$.
\end{proof}

More generally, for a polynomial $f(x) \in \mathbb{Q}[x]$, Bell~\cite{Bell} showed that $\nu_p(f(n))_{n \geq 0}$ is $p$-regular if and only if $f(x)$ factors as a product of linear polynomials in $\mathbb{Q}[x]$, times a polynomial with no root in the $p$-adic integers.

The main purpose of this paper is to prove the following.

\begin{theorem}\label{exact rank}
Let $p \neq 2, 5$ be a prime.
Then $\nu_p(F_{n + 1})_{n \geq 0}$ is a $p$-regular sequence of rank $\alpha(p) + 1$.
\end{theorem}

The sequence $\nu_p(F_{n + 1})_{n \geq 0}$ is $p$-regular for the primes $2$ and $5$ as well.
In Section~\ref{Determining the rank} we prove Theorem~\ref{exact rank} and show that for $p = 2$ the rank is $5$ and for $p = 5$ the rank is $2$.
First, in Section~\ref{Closure properties}, we discuss some closure properties of regular sequences and show that the $p$-regularity of $\nu_p(F_{n + 1})_{n \geq 0}$ follows from these properties, although the upper bound we obtain for the rank is not sharp.

In an unpublished note (with the same title as the current manuscript), we were able to show the $p$-regularity of the $p$-adic valuation of the Fibonacci sequence, although we did not establish the rank.
In \cite{Shu--Yao}, Shu and Yao also proved the $p$-regularity of $\nu_p(F_{n+1})_{n \geq 0}$, using $p$-adic analytic methods.
More generally, they characterized sequences satisfying a linear recurrence of order $2$ with constant coefficients for which the sequence of $p$-adic valuations is $p$-regular.
In particular, not every such sequence of $p$-adic valuations is $p$-regular.

\section{Closure properties}\label{Closure properties}

For a fixed $k$, $k$-regular sequences satisfy several closure properties.
We shall make use of two of these closure properties.
Since we are interested in the rank of $\nu_p(F_{n+1})_{n \geq 0}$, we will record the bound we obtain on the rank of each resulting sequence in terms of the ranks of the initial sequences.
The first property is closure under termwise addition~\cite[Theorem~2.5]{Allouche--Shallit}.

\begin{theorem}\label{addition}
Let $k \geq 2$, and let $s(n)_{n \geq 0}$ and $t(n)_{n \geq 0}$ be $k$-regular sequences.
Then $(s(n) + t(n))_{n \geq 0}$ is a $k$-regular sequence of rank at most $\rank s + \rank t$.
\end{theorem}

\begin{proof}
Let $B$ be a finite set of sequences that generate $\langle \ker_k s(n)_{n \geq 0} \rangle$, and let $C$ be a finite set of sequences that generate $\langle \ker_k t(n)_{n \geq 0} \rangle$.
Then $B \cup C$ generates $\langle \ker_k (s(n) + t(n))_{n \geq 0} \rangle$.
\end{proof}

The second closure property states that riffling together a number of $k$-regular sequences produces a $k$-regular sequence.

\begin{theorem}\label{riffle}
Let $s(n)_{n \geq 0}$ be a sequence.
Let $k \geq 2$ and $a \geq 1$ be relatively prime integers such that $s(a n + b)_{n \geq 0}$ is $k$-regular for each $0 \leq b \leq a - 1$.
Let $s_b(n) = s(a n + b)$.
Then $s(n)_{n \geq 0}$ is $k$-regular of rank at most $a \cdot \sum_{b=0}^{a-1} \rank s_b$.
\end{theorem}

Allouche and Shallit~\cite[Theorem~2.7]{Allouche--Shallit} state this closure property without the condition that $k$ and $a$ are relatively prime, although the proof is incomplete in the case $\gcd(k, a) \geq 2$.

\begin{proof}
For each $0 \leq b \leq a - 1$, let
\[
	t_b(n) =
	\begin{cases}
		s(n)	& \text{if $n \equiv b \mod a$} \\
		0	& \text{if $n \nequiv b \mod a$.}
	\end{cases}
\]
We claim that $t_b(n)_{n \geq 0}$ is $k$-regular for each $b$.
Consider the element $t_b(k^j n + c)_{n \geq 0}$ of the $k$-kernel of $t_b(n)_{n \geq 0}$.
If $k^j n + c \equiv b \mod a$ then $n \equiv k^{-j} (b - c) \mod a$ since $\gcd(k, a) = 1$.
Therefore $t_b(k^j n + c)_{n \geq 0}$ is the sequence
\[
	s\left(k^j \left(a n + \left(k^{-j} (b - c) \bmod a\right)\right) + c\right)_{n \geq 0}
\]
interspersed with blocks of $a - 1$ zeros, beginning at some offset.
In other words, $t_b(k^j n + c)_{n \geq 0}$ is the sequence
\[
	s_b\left(k^j n + \frac{k^j \left(k^{-j} (b - c) \bmod a\right) + c - b}{a}\right)_{n \geq 0}
\]
interspersed with blocks of $a - 1$ zeros, beginning at some offset.
There are $a$ possible offsets.
The module generated by all elements in the $k$-kernel of $t_b(n)_{n \geq 0}$ with a given offset has rank at most $\rank s_b$.
Therefore $t_b(n)_{n \geq 0}$ is $k$-regular with $\rank t_b \leq a \rank s_b$.
Since $s(n) = \sum_{b=0}^{a-1} t_b(n)$, it follows from Theorem~\ref{addition} that $s(n)_{n \geq 0}$ is $k$-regular with $\rank s \leq a \sum_{b=0}^{a-1} \rank s_b$.
\end{proof}

The $p$-regularity of $\nu_p(F_{n+1})_{n \geq 0}$ follows from these closure properties.

\begin{theorem}\label{closure bounds}
Let $p \neq 2, 5$ be a prime.
Then $\nu_p(F_{n + 1})_{n \geq 0}$ is a $p$-regular sequence of rank at most $3 \alpha(p)$.
\end{theorem}

\begin{proof}
By Theorem~\ref{Lengyel} we have
\[
	\nu_p(F_n) =
	\begin{cases}
		\nu_p(n) + \nu_p(F_{\alpha(p)})	& \text{if $n \equiv 0 \mod \alpha(p)$} \\
		0				& \text{if $n \nequiv 0 \mod \alpha(p)$}
	\end{cases}
\]
for $n \geq 1$.
Therefore $\nu_p(F_{n + 1})_{n \geq 0}$ is a riffle of the sequence
\[
	\left(\nu_p(\alpha(p) \cdot (n+1)) + \nu_p(F_{\alpha(p)})\right)_{n \geq 0}
\]
and $\alpha(p) - 1$ zero sequences.
It follows from Theorem~\ref{mod5} that $p$ and $\alpha(p)$ are relatively prime, so $\nu_p(\alpha(p) \cdot (n+1)) = \nu_p(n+1)$.
We have shown in Theorem~\ref{valuation rank} that the rank of $\nu_k(n+1)_{n \geq 0}$ is $2$, and the constant sequence $\nu_p(F_{\alpha(p)})_{n \geq 0}$ has rank $1$, so by Theorem~\ref{addition} their sum has rank at most $3$.
Since the rank of the zero sequence is $0$, Theorem~\ref{riffle} now implies that $\nu_p(F_{n + 1})_{n \geq 0}$ is $p$-regular with rank at most $3 \alpha(p)$.
\end{proof}

\section{Determining the rank}\label{Determining the rank}

In this section we prove Theorem~\ref{exact rank}, showing that the rank of $\nu_p(F_{n + 1})_{n \geq 0}$ is less than the bound $3 \alpha(p)$ given by Theorem~\ref{closure bounds}.
As in the proof of Theorem~\ref{valuation rank}, we exhibit generators and relations for the $\Z$-module generated by the $p$-kernel of $\nu_p(F_{n+1})_{n \geq 0}$.
First, however, we address the primes $2$ and $5$.

\begin{theorem}
The sequence $\nu_5(F_{n + 1})_{n \geq 0}$ is a $5$-regular sequence of rank $2$.
\end{theorem}

\begin{proof}
This follows immediately from Theorem~\ref{Lengyel} and Theorem~\ref{valuation rank}.
\end{proof}

\begin{theorem}
The sequence $\nu_2(F_{n + 1})_{n \geq 0}$ is a $2$-regular sequence of rank $5$.
\end{theorem}

\begin{proof}
Let $B$ be the set
\[
	\left\{
	\nu_2(F_{n + 1})_{n \geq 0}, \,
	\nu_2(F_{2 n + 1})_{n \geq 0}, \,
	\nu_2(F_{2 n + 2})_{n \geq 0}, \,
	\nu_2(F_{4 n + 1})_{n \geq 0}, \,
	\nu_2(F_{4 n + 3})_{n \geq 0}
	\right\}.
\]
We claim $\langle B \rangle = \langle \ker_2 \nu_2(F_{n + 1})_{n \geq 0} \rangle$.
The identities
\begin{align*}
	\nu_2(F_{(2 n + 0) + 1}) &= \nu_2(F_{2 n + 1})		&	\nu_2(F_{(2 n + 1) + 1}) &= \nu_2(F_{2 n + 2}) \\
	\nu_2(F_{2 (2 n + 0) + 1}) &= \nu_2(F_{4 n + 1})	&	\nu_2(F_{2 (2 n + 1) + 1}) &= \nu_2(F_{4 n + 3}) \\
	\nu_2(F_{2 (2 n + 0) + 2}) &= 3 \nu_2(F_{2 n + 1})	&	\nu_2(F_{2 (2 n + 1) + 2}) &= \nu_2(F_{4 n + 1}) + \nu_2(F_{2 n + 2}) \\
	\nu_2(F_{4 (2 n + 0) + 1}) &= \nu_2(F_{2 n + 1})	&	\nu_2(F_{4 (2 n + 1) + 1}) &= \nu_2(F_{4 n + 1}) \\
	\nu_2(F_{4 (2 n + 0) + 3}) &= \nu_2(F_{4 n + 3})	&	\nu_2(F_{4 (2 n + 1) + 3}) &= \nu_2(F_{2 n + 1})
\end{align*}
follow from applications of Theorem~\ref{Lengyel} and show that for each sequence $s(n)_{n \geq 0} \in B$, $s(2 n + i)_{n \geq 0} \in \langle B \rangle$ for each $0 \leq i \leq 1$.
Finally, one checks that the sequences in $B$ are linearly independent (for example, by computing the first $16$ terms of each).
\end{proof}

We now prove Theorem~\ref{exact rank}.
To make use of Theorem~\ref{mod5}, we break the remaining primes into equivalence classes modulo $5$.

\begin{theorem}\label{1,4}
Let $p$ be a prime such that $p \equiv 1,4 \mod 5$.
Then the $p$-regular rank of $\nu_p(F_{n + 1})_{n \geq 0}$ is $\alpha(p) + 1$.
\end{theorem}

\begin{proof}
Recall that $p \equiv 1 \mod \alpha(p)$ by Theorem~\ref{mod5}.
For $n \geq 0$ and $0 \leq j \leq \alpha(p) - 1$, let $s_j(n) = \frac{1}{\nu_p(F_{\alpha(p)})} \nu_p(F_{p n + j + 1})$.
By Theorem~\ref{Lengyel},
\[
	s_j(n) =
	\frac{1}{\nu_p(F_{\alpha(p)})}
	\begin{cases}
		\nu_p(p n + j + 1) + \nu_p(F_{\alpha(p)})	& \text{if $p n + j + 1 \equiv 0 \mod \alpha(p)$} \\
		0						& \text{if $p n + j + 1 \nequiv 0 \mod \alpha(p)$.}
	\end{cases}
\]
Since $p \nmid j + 1$, we have $\nu_p(p n + j + 1) = 0$, and therefore
\[
	s_j(n) =
	\begin{cases}
		1	& \text{if $n \equiv -(j + 1) \mod \alpha(p)$} \\
		0	& \text{if $n \nequiv -(j + 1) \mod \alpha(p)$.}
	\end{cases}
\]
In particular, each sequence $s_j(n)_{n \geq 0}$ is a sequence of integers.
Consider the set $B$ of size $\alpha(p) + 1$ consisting of the original sequence $\nu_p(F_{n+1})_{n \geq 0}$ and the sequences $s_j(n)_{n \geq 0}$ for $0 \leq j \leq \alpha(p) - 1$.

First we show that the $p$ subsequences of $\nu_p(F_{n+1})_{n \geq 0}$ are $\Z$-linear combinations of elements of $B$.
We claim for $n \geq 0$ that
\[
	\nu_p(F_{p n + i + 1}) =
	\begin{cases}
		\nu_p(F_{\alpha(p)}) s_{i \bmod \alpha(p)}(n)	& \text{if $0 \leq i \leq p - 2$} \\
		\nu_p(F_{n + 1}) + s_0(n)			& \text{if $i = p-1$.}
	\end{cases}
\]
To see that this holds for $0 \leq i \leq p - 2$, apply Theorem~\ref{Lengyel} to the left side and use the fact that $\nu_p(p n + i + 1) = 0$ since $p \nmid i + 1$.
For $i = p - 1$, both sides are equal (again by Theorem~\ref{Lengyel}) to
\[
	\begin{cases}
		\nu_p(n+1) + \nu_p(F_{\alpha(p)}) + 1	& \text{if $n+1 \equiv 0 \mod \alpha(p)$} \\
		0					& \text{if $n+1 \nequiv 0 \mod \alpha(p)$.}
	\end{cases}
\]

Next we show that the $p$ subsequences of $s_j(n)_{n \geq 0}$ for each $0 \leq j \leq \alpha(p) - 1$ are $\Z$-linear combinations of elements of $B$.
We claim
\[
	s_j(p n + i) = s_{i + j \bmod \alpha(p)}(n)
\]
for $n \geq 0$ and $0 \leq i \leq p - 1$.
Indeed, both sides are equal to
\[
	\begin{cases}
		1	& \text{if $n + i + j + 1 \equiv 0 \mod \alpha(p)$} \\
		0	& \text{if $n + i + j + 1 \nequiv 0 \mod \alpha(p)$.}
	\end{cases}
\]

We have shown that the rank of $\nu_p(F_{n + 1})_{n \geq 0}$ is at most $\alpha(p) + 1$.
It remains to show that there are $\alpha(p) + 1$ linearly independent sequences in the $p$-kernel of $\nu_p(F_{n + 1})_{n \geq 0}$.
Since each sequence in $B$ is a scalar multiple of a sequence in the $p$-kernel, it suffices to show that $B$ is linearly independent.
Clearly the $s_j(n)_{n \geq 0}$ for $0 \leq j \leq \alpha(p) - 1$ are linearly independent, and the sequence $\nu_p(F_{n+1})_{n \geq 0}$ is not a linear combination of the $s_j(n)_{n \geq 0}$ since $\nu_p(F_{p \alpha(p)}) = 1 + \nu_p(F_{\alpha(p)}) \neq \nu_p(F_{\alpha(p)})$ even though $p \alpha(p) \equiv \alpha(p) \mod \alpha(p)$.
\end{proof}

Note that in the proof of Theorem~\ref{1,4} some effort was required to accommodate the possibility that $\nu_p(F_{\alpha(p)}) > 1$.
If no Wall--Sun--Sun primes exist, then $\nu_p(F_{\alpha(p)}) = 1$ for all $p$ and some of the details could be simplified.
The same is true for the proof of the following theorem.

\begin{theorem}\label{2,3}
Let $p \neq 2$ be a prime such that $p \equiv 2,3 \mod 5$.
Then the $p$-regular rank of $\nu_p(F_{n + 1})_{n \geq 0}$ is $\alpha(p) + 1$.
\end{theorem}

\begin{proof}
By Theorem~\ref{mod5}, we have $p \equiv -1 \mod \alpha(p)$.
We construct a basis as in the proof of Theorem~\ref{1,4}.
However, there is a minor case distinction to be made.
If $\alpha(p) < p + 1$, let $s_j(n) = \frac{1}{\nu_p(F_{\alpha(p)})} \nu_p(F_{p n + j + 1})$ for $0 \leq j \leq \alpha(p) - 1$ as before.
If $\alpha(p) = p + 1$, let $s_j(n) = \frac{1}{\nu_p(F_{\alpha(p)})} \nu_p(F_{p n + j + 1})$ only for $0 \leq j \leq \alpha(p) - 3$; then let $s_{\alpha(p) - 2}(n) = \frac{1}{\nu_p(F_{\alpha(p)})} \nu_p(F_{p^2 n + 1})$ and $s_{\alpha(p) - 1}(n) = \frac{1}{\nu_p(F_{\alpha(p)})} \nu_p(F_{p^2 n + p + 1})$.
Then by Theorem~\ref{Lengyel} (using the fact that $p \neq 2$)
\[
	s_j(n) =
	\begin{cases}
		1	& \text{if $n \equiv j + 1 \mod \alpha(p)$} \\
		0	& \text{if $n \nequiv j + 1 \mod \alpha(p)$}
	\end{cases}
\]
for $0 \leq j \leq \alpha(p) - 1$.

Let the set $B$ consist of the original sequence $\nu_p(F_{n+1})_{n \geq 0}$ and the sequences $s_j(n)_{n \geq 0}$ for $0 \leq j \leq \alpha(p) - 1$.
One checks that the relations
\begin{align*}
	\nu_p(F_{p n + i + 1}) &=
	\begin{cases}
		\nu_p(F_{\alpha(p)}) s_{i \bmod \alpha(p)}(n)	& \text{if $0 \leq i \leq p - 2$} \\
		\nu_p(F_{n + 1}) + s_{\alpha(p) - 2}(n)			& \text{if $i = p-1$,}
	\end{cases} \\
	s_j(p n + i) &= s_{i - j - 2 \bmod \alpha(p)}(n)
\end{align*}
for $0 \leq i \leq p - 1$ and $0 \leq j \leq \alpha(p) - 1$ follow from Theorem~\ref{Lengyel}.
We omit the details.
The proof that $B$ is linearly independent is identical to that in the proof of Theorem~\ref{1,4}.
\end{proof}

We have now proved Theorem~\ref{exact rank}.


\end{document}